\documentclass{amsart}
\usepackage{amssymb}
\usepackage{enumerate}
\def\AA{\mathbb A}
\def\bu{$\bullet$\quad}
\def\CA{\mathcal A}
\def\casesep{,&\text{if }}
\def\CB{\mathcal B}
\def\CC{\mathbb C}
\def\CO{\mathcal O}
\def\DD{\mathbb D}
\def\Delta{\varDelta}

\def\Log{\operatorname{Log}}
\def\Omega{\varOmega}
\def\phi{\varphi}
\def\Phi{\varPhi}

\def\rho{\varrho}
\def\SH{\mathcal{SH}}
\def\too{\longrightarrow}
\def\tuu{\longmapsto}
\def\TT{\mathbb T}
\def\wdtl{\widetilde}
\def\wdht{\widehat}
\makeatletter
\def\th@mytheorem{%
  \let\thm@indent\noindent
  \thm@headfont{\bfseries}
    \itshape
}
\def\th@myremark{%
  \let\thm@indent\noindent
  \thm@headfont{\bfseries}
}
\makeatother
\theoremstyle{mytheorem}
\newtheorem{Theorem}{Theorem}
\theoremstyle{myremark}
\newtheorem{Remark}[Theorem]{Remark}
\begin{document}
\title[A counterexample --- revisited]{A counterexample to a theorem of Bremermann on Shilov boundaries --- revisited}

\author[M.~Jarnicki]{Marek Jarnicki}
\address{Jagiellonian University, Faculty of Mathematics and Computer Science, Institute of Mathematics,
{\L}ojasiewicza 6, 30-348 Krak\'ow, Poland}
\email{Marek.Jarnicki@im.uj.edu.pl}

\author[P.~Pflug]{Peter Pflug}
\address{Carl von Ossietzky Universit\"at Oldenburg, Institut f\"ur Mathematik,
Postfach 2503, D-26111 Oldenburg, Germany}
\email{Peter.Pflug@uni-oldenburg.de}

\thanks{The research was partially supported by grant no.~UMO-2011/03/B/ST1/04758 of the Polish National Science Center (NCN)}

\begin{abstract}
We continue to discuss the example presented in \cite{JarPfl2015}. In particular, we clarify some gaps and complete the description of the Shilov boundary.
\end{abstract}



\maketitle

For a bounded domain $D\subset\CC^n$ let $\CA(D)$ (resp.~$\CO(\overline D)$) denote the space of all continuous functions $f:\overline D\too\CC$ such that $f|_D$ is holomorphic
(resp.~$f$ extends holomorphically to a neighborhood of $\overline D$). Let $\partial_SD$ (resp.~$\partial_BD$) be the \emph{Shilov} (resp.~\emph{Bergman})
\emph{boundary} of $D$, i.e.~the minimal compact set $K\subset\overline D$ such that $\max\limits_K|f|=\max\limits_{\overline D}|f|$ for every $f\in\CA(D)$
(resp.~$f\in\CO(\overline D)$). Obviously, $\CO(\overline D)\subset\CA(D)$ and hence $\partial_BD\subset\partial_SD\subset\partial D$.
Notice that, in general, $\partial_BD\varsubsetneq\partial_SD$, e.g.~for the domain $D:=\{(z,w)\in\CC^2: 0<|z|<1,\;|w|<|z|^{-\log|z|}\}$ (cf.~\cite{Fuk1965}, \S\;16).

The algebra $\CA(D)$ (resp.~$\CB(D):=$ the uniform closure in $\CA(D)$ of $\CO(\overline D)$) endowed with the supremum norm is a Banach algebra.
Then $\partial_SD$ (resp.~$\partial_BD$) coincides with the Shilov boundary of $\CA(D)$ (resp.~$\CB(D)$ in the sense of uniform
algebras (cf.~\cite{Bis1959}). Note that the peak points of $\CA(D)$ (resp.~$\CB(D))$ are dense in $\partial_SD$
(resp.~$\partial_B(D)$)) (cf.~\cite{Bis1959}). Recall that a point $a\in\overline D$ is called a peak point for $\CA(D)$ (resp.~$\CB(D)$)
if there is an $f\in\CA(D)$ (resp.~$\CB(D)$) with $f(a)=1$ and $|f(z)|<1$ for all $z\in\overline D\setminus\{a\}$; $f$ is called an associated \emph{peak function}.

Assume that the envelope of holomorphy $\wdtl D$ of $D$ is univalent. In \cite{JarPfl2015} we were interested in
answering whether $\partial_SD=\partial_S\wdtl D$ (resp.~$\partial_BD=\partial_B\wdtl D$).

\begin{Remark}
Notice that:
\begin{itemize}
\item $\partial_S\wdtl D\subset\partial_SD$,

\item $\partial_B\wdtl D\subset\partial_BD$,

\item if $\CA(D)\subset\CA(\wdtl D)|_{\overline D}$ (resp.~$\CO(\overline D)\subset\CO(\overline{\wdtl D})|_{\overline D}$), then $\partial_SD=\partial_S\wdtl D$
(resp.~$\partial_BD=\partial_B\wdtl D$).
\end{itemize}
\end{Remark}

In \cite{JarPfl2015} we studied the following  bounded Hartogs domain $D\subset\CC^2$:
\begin{multline*}
D:=\{(re^{i\phi},w)\in\CC^2: \tfrac12<r<1,\;\phi\in(0,2\pi),\\
\begin{cases}
0<\phi\leq\frac\pi2 &\!\!\!\!\Longrightarrow |w|<1\\
\frac\pi2<\phi<\frac{3\pi}2 &\!\!\!\!\Longrightarrow |w|<3\\
\frac{3\pi}2\leq\phi<2\pi &\!\!\!\!\Longrightarrow 2<|w|<3\\
\end{cases}\};
\end{multline*}
it is known that $D$ has a univalent envelope of holomorphy $\wdtl D$. The main result of \cite{JarPfl2015} is the following theorem.
\begin{Theorem}
$\partial_S\wdtl D\varsubsetneq\partial_SD$,  $\partial_B\wdtl D\varsubsetneq\partial_BD$, and
$\CO(\overline D)\setminus\CA(\wdtl D)|_{\overline D}\neq\varnothing$.
\end{Theorem}

The proof consists of the following two parts:
\begin{enumerate}
\item $\partial_S\wdtl D\cap(I\times\DD(3))=\varnothing$, where $I:=[\frac12,1]$, $\DD$ is the unit disc, and $\DD(r):=r\DD$.

\item There exists a function $h\in\CO(\overline D)$ (effectively given) such that
$$
h(x,w)=\begin{cases}
e^{-2\pi+i\log x}\casesep (x,w)\in I\times\overline\AA(2,3)\\
e^{i\log x}\casesep (x,w)\in I\times\overline\DD
\end{cases}
$$
and $|h|<1$ on the remaining part of $\overline D$, where $\AA(r_-,r_+):=\DD(r_+)\setminus\overline\DD(r_-)$.
\end{enumerate}

Unfortunately, the proof of (1) contains a gap.
The aim of the present note is to close the above gap and to prove some new results related to the Shilov and Bergman boundaries of $D$ and $\wdtl D$.

\medskip

Let
$$
A:=\{z\in\CC: \tfrac12<|z|<1\},\quad I_0:=(\tfrac12,1),\quad A_0:=A\setminus I_0.
$$
By the Cauchy integral formula each function $f\in\CA(D)$ extends holomorphically to the domain
$$
G=\{(z,w)\in A_0\times\CC: |w|e^{V(z)}<1\},
$$
where
$$
V(re^{i\phi}):=\begin{cases}0\casesep 0<\phi\leq\frac\pi2\\-\log 3\casesep \frac\pi2<\phi<2\pi\end{cases}.
$$
Hence (by \cite{JarPfl2000}, Corollary 3.2.18) the envelope of holomorphy $\wdtl D$ is univalent and
$$
\wdtl D=\wdtl G=\{(z,w)\in A_0\times\CC: |w|e^{\wdtl V(z)}<1\},
$$
where
$$
\wdtl V(z):=\sup\{u\in\SH(A_0): u\leq V\}.
$$
Notice that, by the maximum principle for subharmonic functions, we have  $\wdtl V(z)<0$, $z\in A_0$.
Thus, $\partial D\cap(U\times\DD(3))\subset\wdtl D$, where $U:=\{re^{i\phi}: r\in I_0,\;0<\phi<\frac\pi2\}$.
Hence $\partial\wdtl D\cap(U\times\DD(3))$  does not contain points of $\partial_S\wdtl D$.

We are going to prove the following theorem.

\begin{Theorem}\label{ThmMain}
\begin{enumerate}[{\rm(a)}]
\item\label{ThmMaina} $\partial_S\wdtl D\cap(I_0\times\DD(3))=\varnothing$.
\item\label{ThmMainb} For any $a\in I$ there exists a $g=g_a\in\CO(\overline D)$ such that $g(a,w)=1$ for all $w\in\overline\DD$, and $|g|<1$ on
$\overline D\setminus(\{a\}\times\overline\DD)$. In particular, $\partial_BD\cap(\{a\}\times\overline\DD)\neq\varnothing$.
\item\label{ThmMainc} $\partial_BD\supset \{a\}\times\TT$ for every $a\in I_0$, where $\TT:=\partial\DD$. Therefore, $\partial_BD\supset I\times\TT$.
\item\label{ThmMaind} $\displaystyle{\partial_BD\setminus((iI_0)\times(3\TT))=\partial_SD\setminus((iI_0)\times(3\TT))}$
\begin{align*}
&\hskip40pt=\{re^{i\phi}: r\in\{\tfrac12,1\},\;\tfrac\pi2\leq\phi\leq2\pi\}\times(3\TT)\\
&\hskip40pt\cup I_0\times(3\TT)\\
&\hskip40pt\cup I_0\times\TT\\
&\hskip40pt\cup \{re^{i\phi}: r\in\{\tfrac12,1\},\;0\leq\phi\leq\tfrac\pi2\}\times\TT\\
&\hskip40pt=:M_1\cup M_2\cup M_3\cup M_4.
\end{align*}
\end{enumerate}
\end{Theorem}

\begin{Remark}
\begin{enumerate}[(i)]
\item Observe that \eqref{ThmMaina} and \eqref{ThmMainb} close gaps in our former proof.

Indeed, we get $\varnothing\neq\partial_BD\setminus\partial_S\wdtl D\subset
(\partial_SD\setminus\partial_S\wdtl D)\cap(\partial_BD\setminus\partial_B\wdtl D)$. Hence $\partial_S\wdtl D\varsubsetneq\partial_SD$ and
$\partial_B\wdtl D\varsubsetneq\partial_BD$.

Moreover, if $a\in I_0$, then $g_a\in\CO(\overline D)\setminus\CA(\wdtl D)|_{\overline D}$.
\item It seems to be an \emph{open problem} whether $(iI_0)\times(3\TT)\subset\partial_SD$ (resp.~$(iI_0)\times(3\TT)\subset\partial_BD$).
\end{enumerate}
\end{Remark}

\begin{proof}[{\it Proof of Theorem \ref{ThmMain}}]
First, let us make the following elementary observation.

(*) Let $\Sigma$ be an open subset of the boundary $\partial\Omega$ of a bounded domain $\Omega\subset\CC^n$. Suppose that
$\max_{\partial\Omega}|f|=\max_{\partial\Omega\setminus\Sigma}|f|$ for every $f\in\CA(\Omega)$ (resp.~$\CO(\overline\Omega)$). Then
$\partial_S\Omega\cap\Sigma=\varnothing$ (resp.~$\partial_B\Omega\cap\Sigma=\varnothing$).

\medskip

\eqref{ThmMaina} For all $a\in I_0$ and $f\in\CA(\wdtl D)$ the function $f(a,\cdot)$ extends holomorphically to $\DD(3)$.

Indeed, we may define
$\wdht f(z,w):=\frac1{2\pi}\int_{|\zeta|=5/2}\frac{f(z,\zeta)}{\zeta-w}d\zeta$, $z\in\AA(\frac12,1), \frac\pi2<\arg z\leq 2\pi$, $|w|<\frac52$.
Then $\wdht f$ is holomorphic and coincide with $f$ when $\frac\pi2<\arg z<\pi$. Hence using identity theorem we see that $f=\wdht f$ on their common
domain of definition. Using continuity of $f$ we get the claimed extension of $f(a,\cdot)$.

In particular,
$\max_{\{a\}\times\overline\AA(1,3)}|f(a,\cdot)|=\max_{\{a\}\times3\TT}|f(a,\cdot)|$. Hence, by (*) with $\Omega:=\wdtl D$ and $\Sigma:=I_0\times\AA(1,3)$,
we conclude that $\partial_S\wdtl D\cap(I_0\times\AA(1,3))=\varnothing$.
The same argument shows that $\partial_S\wdtl D\cap (I_0\times\DD)=\varnothing$.
(Note that these first two cases can be also handled using the density of the peak points --- see the argument in the next case.)

Suppose that $(z_0,w_0)\in(I_0\times\TT)\cap\partial_S\wdtl D$. Then there is a peak point $(z_1,w_1)$ nearby. Let $f\in\CA(\wdtl D)$ be a function peaking there.
The maximum principle excludes the situation where $z_1\in I_0$. Thus $z_1\in U$, but we already know that $\partial_S\wdtl D\cap(U\times\DD(3))=\varnothing$,
so it is impossible.

Finally, $\partial_S\wdtl D\cap(I_0\times\DD(3))=\varnothing$.

\medskip

\eqref{ThmMainb} Fix an $a\in I$ and let $h$ be as in (2), $w_0:=e^{i\log a}\in\TT$. Define $\phi(w):=\frac12(w+w_0)$, $g:=\phi\circ h$. It is obvious that
$g\in\CO(\overline D)$, $g(a,w)=1$ for all $w\in\overline\DD$, and $|g|<1$ on $\overline D\setminus(\{a\}\times\overline\DD)$.

\medskip

\eqref{ThmMainc} Using \eqref{ThmMaina} we have $\partial_BD\cap(I_0\times\DD)=\varnothing$. Hence, by \eqref{ThmMainb},
$\partial_BD\cap(\{a\}\times\TT)\neq\varnothing$ for every $a\in I_0$. Now using rotational invariance in the second variable of $\partial_BD$
leads to $\{a\}\times\TT\subset\partial_BD$ for all $a\in I_0$.

\medskip

\eqref{ThmMaind} Notice that also $\partial_SD$ is invariant under rotations of the second variable.

\bu Every point from $\partial A\times3\TT$ is a peak point for $\CO(\overline{A\times\DD(3)})$.

Indeed, fix a point $(a,b)\in\partial A\times 3\TT$. Then $a$ is a peak point for $\CO(\overline A)$ and $b$ is a peak point for
$\CO(\overline{\DD(3)})$. So it suffices to take the product of the corresponding peak functions to see that $(a,b)$ is a peak point
for $\CO(\overline{A\times\DD(3)})$.

Thus $M_1\subset\partial_BD\subset\partial_SD$.

\bu Consider the holomorphic function
$$
\overline D\ni (z,w)\overset{\Phi}\tuu\Log(z)\in R:=[-\log2,0]\times[0,2\pi],
$$
where $\Log$ is a branch of logarithm with $\Log(-1)=\pi$.
Note that $\Phi\in\CO(\overline D)$.  For every $a\in I_0$ we have
$\Phi(a,w)\in(-\log2,0)\times\{2\pi\}\in\partial R$ whenever $w\in\overline\AA(2,3)$. It is clear that there exists a function
$\psi_a\in\CO(\overline R)$ such that $\psi_a(\Phi(a,w))=1$, $w\in\overline\AA(2,3)$, and $|\psi_a|<1$ on $\overline R\setminus\{\Phi(a,w)\}$.
Then the function
$$
\overline D\ni(z,w)\overset{f_a}\tuu\psi_a(\Phi(z,w))
$$
may be considered as a function of class $\CO(\overline D)$.
Observe that $f_a(a,w)=1$ for all $w\in\overline\AA(2,3)$,
and $|f_a|<1$ on $\overline D\setminus(\{a\}\times\overline\AA(2,3))$. Fix a $b\in3\TT$ and define $g(z,w):=f_a(z,w)\frac{1+w/b}2$.
Then $g\in\CO(\overline D)$ and $g$ peaks at $(a,b)$. Consequently, $M_2\subset\partial_BD\subset\partial_SD$.

\bu By \eqref{ThmMainc}, $M_3\subset\partial_BD\subset\partial_SD$.

\bu For every $a=re^{i\phi}$ with $r\in\{\frac12,1\}$, $0<\phi<\frac\pi2$ there exists a function $\psi\in\CO(\overline A)$ such that $\psi(a)=1$ and $|\psi|<1$ on
$\overline A\setminus\{a\}$. Hence $\partial_BD\cap(\{a\}\times\overline\DD)\neq\varnothing$. Now, by (*) with $\Omega:=D$,
$\Sigma:=\{re^{i\phi}: r\in\{\frac12,1\},\; 0<\phi<\frac\pi2\}\times\DD$, we conclude that $M_4\subset\partial_BD\subset\partial_SD$.

\medskip

The remaining part of $\partial D$, i.e.~the set $\Sigma:=\partial D\setminus(M_1\cup M_2\cup M_3\cup M_4\cup ((iI_0)\times(3\TT)))$,  is open in $\partial D$.
It remains to use (*).
\end{proof}

\begin{Remark}
We will try to complete the description of $\partial_BD$, $\partial_SD$, $\partial_B\wdtl D$, and $\partial_S\wdtl D$. Any help by the reader will be appreciated.
\end{Remark}

\bibliographystyle{amsplain}

\end{document}